# MODELING CREDIT RISK WITH PARTIAL INFORMATION

By Umut Cetin, Robert Jarrow, Philip Protter and Yildiray Yildirim

*Cornell University, Cornell University, Cornell University and Syracuse University*

This paper provides an alternative approach to Duffie and Lando [*Econometrica* **69** (2001) 633–664] for obtaining a reduced form credit risk model from a structural model. Duffie and Lando obtain a reduced form model by constructing an economy where the market sees the manager's information set plus noise. The noise makes default a surprise to the market. In contrast, we obtain a reduced form model by constructing an economy where the market sees a reduction of the manager's information set. The reduced information makes default a surprise to the market. We provide an explicit formula for the default intensity based on an Azéma martingale, and we use excursion theory of Brownian motions to price risky debt.

**1. Introduction.** Reduced form models have become important tools in the risk management of credit risk [for background references see Jarrow and Yu (2001) and Bielecki and Rutkowski (2002)]. One reason for this is that they usually provide a better fit to market data than structural models do [see Jones, Mason and Rosenfeld (1984), Jarrow, van Deventer and Wang (2002) and Eom, Helwege and Huang (2000)]. Reduced form models take a firm's default process as exogenous with the time of default a stopping time. When in addition the time is totally inaccessible, the market cannot predict the time of default. Yet, managers working within a firm surely know when default is imminent. From a manager's perspective, default is an accessible stopping time (predictable). Usually, in the structural approach default occurs when the firm's value, a continuous sample path process, hits a barrier. This formulation is consistent with the manager's perspective but inconsistent with reduced form models.

Duffie and Lando (2001) link the two perspectives by introducing noise into the market's information set, transforming the manager's accessible default time from the structural approach into the market's inaccessible default time of a reduced form model. Duffie and Lando postulate that the market can only observe the firm's asset value plus noise at equally spaced, discrete time points (and not continuously). And, when default occurs, the market is immediately informed. This noise generates the market's surprise with respect to default, because the firm could nearly be in default (just about to hit the barrier) and the market not yet aware of its imminence. Kusuoka (1999) extends Duffie and Lando's model to continuous time observations of the firm's asset value plus noise. Kusuoka's solution is an application of continuous time filtering theory.

This approach to constructing a reduced form credit model presumes that the market has the same information set as the firm's management, but with noise appended. (Filtering theory was originally formulated for electronic signal processing where the physical problem corresponds to a situation where an electronic signal is received with noise and the noise needs to be "filtered" out.) An interpretation is that accounting reports and/or management press releases either purposefully (e.g., Enron) or inadvertently add extraneous information that obscures the market's knowledge of the firm's asset value. Management knows the firm's value (because this knowledge determines default), but they cannot (or will not) make it known to the market. The market's task is to remove this extraneous noise. Although possible in many situations, this characterization of management's information versus the market's is not exhaustive. An alternative and equally plausible characterization is that the market has the same information as a firm's management, but just less of it. Accounting reports and/or management press releases provide just a reduced set of the information that is available.

Consistent with this alternative perspective, we provide a second approach to the construction of a reduced form credit risk model from a structural model. In our approach, the firm's cash flows, a continuous sample path process, provide the sufficient statistic for default. If the firm's cash flows remain negative for an extended period of time, the firm after exhausting both its lines of credit and easily liquidated assets, defaults. Management observes the firm's cash flows. In contrast, the market observes only a very coarse partitioning of the manager's information set. The market knows only that the cash flow is negative; the firm is experiencing financial distress and the duration of the negative cash flow event, nothing else. This information structure has default being an accessible stopping time for management, but an inaccessible stopping time for the market, yielding the reduced form credit risk model.

To illustrate the economic concepts involved, this paper concentrates on developing a specific example to obtain analytic results. The analytic results solidify intuition and make the economic arguments more transparent. Generalizations and extensions will be readily apparent once the example is well understood. It is our hope that this paper will motivate additional research into this area. Our example provides an explicit representation of the firm's default intensity using an Azéma's martingale [see Emery (1989)]. To illustrate the usefulness of this result, we compute the value of a risky zero-coupon bond using excursion theory of Brownian motions. For another application of excursion theory to option pricing see Chesney, Jeanblanc-Picqué and Yor (1997).

An outline for this paper is as follows. Section 2 presents the structural model. Section 3 presents the reduced form model, Section 4 values a risky zero-coupon bond in the reduced form model, while Section 5 concludes the paper.

**2. The structural model.** We consider a continuous trading economy with a money market account where default-free zero-coupon bonds are traded. In this economy there is a risky firm with debt outstanding in the form of zero-coupon bonds. The details of these traded assets are not needed now, but will be provided later as necessity dictates. The market for these traded securities is assumed to be arbitrage free, but not necessarily complete.

We begin with a filtered probability space $(\Omega, \mathcal{F}, (\mathcal{F}_t)_{0 \leq t \leq T}, \mathbb{Q})$ satisfying the usual conditions. Time $T > 0$ is the final date in the model. The probability $\mathbb{Q}$ is an equivalent martingale probability measure under which the normalized prices of the traded securities follow a martingale. Normalization is by the value of the money market account. The no-arbitrage assumption guarantees the existence, but not the uniqueness of such a probability measure [see Duffie (1996)].

2.1. *Management's information.* Let $X$ be the cash balances of the firm, normalized by the value of the money market account, with the following stochastic differential equation:

$$(2.1) \qquad dX_t = \sigma \, dW_t, \qquad X_0 = x$$

with $x > 0$, $\sigma > 0$, and where $W$ is a standard Brownian motion on the given probability space.

The cash balances of the firm are initialized at $x > 0$ units of the money market account. One should interpret this quantity as the "target" or "optimal" cash balances for the firm. An optimal cash balance could exist because if the firm holds too much cash, it forgoes attractive investment projects and incurs increased tax liabilities, while if it has too little cash, it increases the likelihood of bankruptcy and the occurrence of third party costs [see Brealey and Myers (2001) for related discussion]. The firm attempts to maintain cash balances at this target level, but fluctuations occur due to its operating needs, for example, meeting payrolls, paying suppliers, receiving payments from accounts receivable, and so on. However, without loss of generality, to simplify the presentation we assume that $x = 0$ and $\sigma = 1$, as well.

Under the martingale measure, cash balances have no drift term. Under the empirical measure, however, one would expect that the cash balances should drift at the spot rate of interest. This is consistent with the firm holding its cash balances in the money market account and trying to maintain the target level balance.

The firm's management observes the firm's cash balances. Cash balances can be positive, zero or negative. Negative cash balances correspond to situations where payments owed are not paid, and the firm is in financial distress.

2.2. *The default process.* Let $\mathcal{Z} := \{t \in [0, T] : X(t) = 0\}$ denote the times when the firm's cash balances hit zero. When the cash balances hit zero, the firm has no cash left for making current payments owed. The firm is in financial distress. With zero or negative cash balances, debt payments can only be made by liquidating the firm's assets or by accessing bank lines of credit. The firm can exist with negative cash balances for only a limited period of time. We now formalize this default process.

Associated with the zero set, we define the following function:

$$g(t) := \sup\{s \le t : X_s = 0\}.$$

The random time $g(t)$ corresponds to the last time (before $t$) that cash balances hit zero. Let

$$\tau_\alpha := \inf\left\{t > 0 : t - g(t) \ge \frac{\alpha^2}{2}, \text{ where } X_s < 0 \text{ for } s \in (g(t-), t)\right\}$$

for some $\alpha \in \mathbb{R}_+$ be the random time that measures the onset of a potential default situation for the firm. Formally, $\tau_\alpha$ is the first time that the firm's cash balances have continued to be negative for at least $\alpha^2/2$ units of time. The constant $\alpha$ is a parameter of the default process (that could be estimated from market data). We let $\tau$ denote the time of default. We assume that

$$\tau := \inf\{t > \tau_\alpha : X_t = 2X_{\tau_\alpha}\}.$$

Default occurs the first time, after $\tau_\alpha$, that the cash balances double in magnitude. The intuition is that after being below zero for a long time, the firm uses up all its slack (lines of credit, etc.) to meet its debt payments. If it ever hits $2X_{\tau_\alpha}$ afterwards, it has no slack left, so it defaults. The doubling in absolute magnitude of the cash balances prior to default is only for analytic convenience, and it has no economic content. The generalization of this assumption is a subject for future research. The above process is what the firm's management observes.

**3. The reduced form model.** This section studies the structural model under the market's information set. It is shown here that the bankruptcy process, as viewed by the market, follows a reduced form model where the indicator function of the default time is a point process with an intensity.

In contrast to the manager's information, the market does not see the firm's cash balances. Instead, until the firm has had prolonged negative cash balances for a certain time, that is, until random time $\tau_\alpha$, the market only knows when the firm has positive cash balances or when it has negative or zero cash balances, and whether the cash balances are above or below the default threshold $2X_{\tau_\alpha}$ afterwards. In this respect, we introduce a new process:

$$Y_t = \begin{cases} X_t, & \text{for } t < \tau_\alpha, \\ 2X_{\tau_\alpha} - X_t, & \text{for } t \geq \tau_\alpha. \end{cases}$$

Note that $Y$ is also an $\mathcal{F}$-Brownian motion and

$$\tau = \inf\{t \geq \tau_\alpha : Y_t = 0\}.$$

Let

$$\text{sign}(x) = \begin{cases} 1, & \text{if } x > 0, \\ -1, & \text{if } x \leq 0. \end{cases}$$

Set $\tilde{\mathcal{G}}_t := \sigma\{\text{sign}(Y_s); s \leq t\}$ and let $(\mathcal{G}_t)_{0 \leq t \leq T}$ denote the $\mathbb{Q}$-complete and right continuous version of the filtration $(\tilde{\mathcal{G}}_t)_{0 \leq t \leq T}$; $(\mathcal{G}_t)_{0 \leq t \leq T}$ is the information set that the market observes. As seen, the market's information set is a very coarse filtering of the manager's information set. In essence, the market observes when the firm is in financial distress, and the duration of this situation.

Given this information, the market values the firm's liabilities by taking conditional expectations under the martingale measure $\mathbb{Q}$. This valuation is studied in the next section.

We now derive the intensity for the default time as seen by the market. Let $\tilde{Y}_t = 2/\sqrt{\pi} Y_t$. Then signs and zero sets of $Y$ and $\tilde{Y}$ are the same. Define $M_t := E[\tilde{Y}_t | \mathcal{G}_t]$. Then, $M$ is the Azéma's martingale on $(\Omega, (\mathcal{G}_t)_{0 \leq t \leq T}, \mathbb{Q})$. [Note that Azéma's martingale has already been used in finance, but in a different context; see Dritschel and Protter (1999).] Its quadratic variation satisfies the following "structure equation":

$$(3.1) \qquad d[M, M]_t = dt - M_{t-}\, dM_t.$$

Azéma's martingale is a strong Markov process. For an extensive treatment of Azéma's martingale and the structure equation, see Emery (1989). We also have the following formula for $M$:

$$(3.2) \qquad M_t = \text{sign}(Y_t)\sqrt{2}\sqrt{t - \bar{g}_t},$$

where $\bar{g}_t := \sup\{s \leq t : Y_s = 0\}$. It is easily seen that $\tau$ can be equivalently written as

$$\tau = \inf\{t > 0 : \Delta M_t \geq \alpha\}.$$

Therefore, $\tau$ is a jump time of Azéma's martingale, hence it is totally inaccessible in the filtration $(\mathcal{G}_t)_{0 \leq t \leq T}$. Also note that $\tau_\alpha = \inf\{t > 0 : M_{t-} \leq -\alpha\}$.

However, $\Delta M_t = -M_{t-}\mathbb{1}_{[M_{t-} \neq M_t]}$. So, $\tau_\alpha \leq \tau$ a.s. Furthermore, $\tau_\alpha$ is a predictable stopping time which implies $\mathbb{Q}[\tau = \tau_\alpha] = 0$. Hence, $\tau_\alpha < \tau$ a.s.

Define $N_t := \mathbb{1}_{[t \geq \tau]}$. By the Doob–Meyer decomposition [see, e.g., Protter (1990), page 90], there exists a continuous, increasing, and predictable (also known as locally natural) process, $A$, such that $N - A$ is a $\mathcal{G}$-martingale which has only one jump, at $\tau$, and of size equal to 1.

THEOREM 3.1. $\tau$ has a $\mathcal{G}$-intensity, that is, $A$ is of the form $A_t = \int_0^{t \wedge \tau} \lambda_s \, ds$. Furthermore, $\lambda_t = \mathbb{1}_{[t > \tau_\alpha]} 1/(2[t - \bar{g}_{t-}])$ for $0 \leq t \leq \tau$, and $\lambda_t = 0$ for $t > \tau$.

PROOF. Let $A_t = A_{t \wedge \tau}$. Then

$$(3.3) \qquad H_t := N_t - A_t = \int_0^{t \wedge \tau} h_s \, dM_s$$

for some $\mathcal{G}$-predictable process $h_s$ since $M$ possesses the predictable representation property. [This is proved in Emery (1989).] Since $A$ is continuous and of finite variation, and $N$ is a quadratic pure jump semimartingale, we have

$$(3.4) \qquad [H, H]_t = [N, N]_t = N_t.$$

Also,

$$
\begin{aligned}
(3.5) \qquad [H, H]_t &= \int_0^{t \wedge \tau} h_s^2 \, d[M, M]_s \\
&= \int_0^{t \wedge \tau} h_s^2 \, ds - \int_0^{t \wedge \tau} h_s^2 M_{s-} \, dM_s,
\end{aligned}
$$

where the second equality follows from (3.1). Combining (3.3)–(3.5) yields

$$\int_0^{t \wedge \tau} h_s^2 \, ds - \int_0^{t \wedge \tau} h_s^2 M_{s-} \, dM_s - A_t = \int_0^{t \wedge \tau} h_s \, dM_s,$$

which implies

$$(3.6) \qquad \int_0^{t \wedge \tau} h_s^2 \, ds - A_t = \int_0^{t \wedge \tau} h_s^2 M_{s-} \, dM_s + \int_0^{t \wedge \tau} h_s \, dM_s.$$

The left-hand side of the previous expression is continuous. Hence

$$(3.7) \qquad \int_0^{t \wedge \tau} (h_s^2 M_{s-} + h_s) \, dM_s = 0.$$

We compute the predictable quadratic variation to get

$$(3.8) \qquad \int_0^{t \wedge \tau} h_s^2 (h_s M_{s-} + 1)^2 \, ds = 0.$$

The optional sampling theorem implies that $(N - A)_{t \wedge \tau_\alpha} = 0$, since $N = 0$ before and at $\tau_\alpha$. Therefore we get $h = 0$ on $[0, \tau_\alpha]$. On the other hand, (3.8) gives $h_s = 0$ or $h_s = -1/M_{s-}$ on $[\tau_\alpha, \tau]$. But (3.3) implies $h_s$ cannot be identically 0 on $(\tau_\alpha, \tau]$, and we see that $h_s = -\mathbb{1}_{[s > \tau_\alpha]} 1/M_{s-}$ satisfies (3.7). Therefore we deduce a version of $H$ which is given by

$$H_t = -\int_0^{t \wedge \tau} \mathbb{1}_{[s > \tau_\alpha]} \frac{1}{M_{s-}} \, dM_s$$

and thus $H$ jumps only at $\tau$ and its jump size is given by

$$\Delta H_\tau = -\mathbb{1}_{[\tau > \tau_\alpha]} \frac{1}{M_{\tau-}} \Delta M_\tau = -\frac{1}{M_{\tau-}} (-M_{\tau-}) = 1.$$

Therefore, (3.6) and (3.7) together imply

$$A_t = \int_0^{t \wedge \tau} h_s^2 \, ds$$
$$= \int_0^{t \wedge \tau} \mathbb{1}_{[s \geq \tau_\alpha]} \frac{1}{M_{s-}^2} \, ds. \qquad \square$$

We choose the intensity equal to 0 after time $\tau$, although other choices might be possible since we are dealing with $\int_0^{t \wedge \tau} \lambda_s \, ds$. This theorem shows that under the market's information set, default is given by a totally inaccessible stopping time, generating a reduced form model from the market's perspective. We have an explicit representation of the intensity process as given by $\lambda_t = \mathbb{1}_{[t > \tau_\alpha]} 1/(2[t - \bar{g}_{t-}])$. The firm's default intensity is zero until time $\tau_\alpha$ is reached. After time $\tau_\alpha$, the default intensity declines with the length of time that the firm remains in financial distress $(t - \bar{g}_{t-})$. The interpretation is that the longer the firm survives in the state of financial distress, the less likely it is to default. Presumably, the firm is more likely to recover and not reach the default magnitude of cash balances given by $2X_{\tau_\alpha}$. With this intensity, the market can value risky bonds and credit derivatives. This valuation is discussed in the next section.

**4. Valuation of a risky zero-coupon bond.** Perhaps one of the most important uses of reduced form credit risk models is to price risky bonds and credit derivatives. This section studies the pricing of risky zero-coupon bonds. Let $(S_t)_{t \in [0,T]}$ denote the price process of a risky zero coupon bond issued by this firm that pays \$1 at time $T$ if no default occurs prior to that date, and zero dollars otherwise. Then, under the no arbitrage assumption, $S$ is given by

$$(4.1) \qquad S_t = \mathbb{E}\Big[\exp\Big(-\int_t^T r_u\,du\Big)\mathbb{1}_{[\tau>T]}\Big|\mathcal{G}_t\Big]\mathbb{1}_{[\tau>t]},$$

where $r_u$ is the instantaneous interest rate at time $u$, and $\mathbb{E}$ refers to the expectation under risk neutral probability law.

To facilitate the evaluation of expression (4.1), we will assume that interest rates are deterministic. In this case, the price of the risky bond becomes

$$(4.2) \qquad S_t = \exp\Big(-\int_t^T r_u\,du\Big)\mathbb{E}[\mathbb{1}_{[\tau>T]}|\mathcal{G}_t]\mathbb{1}_{[\tau>t]}.$$

Let $V_t = \mathbb{1}_{[t<T]}\mathbb{E}[\exp(-\int_t^T \lambda_u\,du)|\mathcal{G}_t]$, where $\lambda$ is the intensity process as given in Theorem 3.1. Duffie, Schroder and Skiadas (1996) give the following formula for $\mathbb{E}[\mathbb{1}_{[\tau>T]}|\mathcal{G}_t]$:

$$\mathbb{E}[\mathbb{1}_{[\tau>T]}|\mathcal{G}_t] = V_t - \mathbb{E}[\Delta V_\tau|\mathcal{G}_t] \qquad \text{on } [t<\tau].$$

The rest of this section is devoted to the computation of this conditional expectation. Define $L_\alpha := \tau - \bar{g}_{\tau_\alpha}$; $L_\alpha$ is the length of the first excursion of Brownian motion below zero exceeding length $\alpha^2/2$. Note that $V_t = 1$ on $[t<T]$ after $\tau$. On $[t<\tau]$,

$$V_t = \mathbb{1}_{[t<\tau]}\mathbb{1}_{[t\geq\tau_\alpha]}\mathbb{E}\Big[\exp\Big(-\int_t^T \lambda_u\,du\Big)\Big|\mathcal{G}_t\Big]$$

$$\qquad + \mathbb{1}_{[t<\tau_\alpha]}\mathbb{E}\Big[\exp\Big(-\int_t^T \lambda_u\,du\Big)\Big|\mathcal{G}_t\Big]$$

$$(4.3) \qquad = +\mathbb{1}_{[t<\tau]}\mathbb{1}_{[t\geq\tau_\alpha]}\sqrt{t-\bar{g}_t}\,\mathbb{E}\Big[\frac{1}{\sqrt{\tau-\bar{g}_t}}\mathbb{1}_{[\tau\leq T]}\Big|\mathcal{G}_t\Big]$$

$$(4.4) \qquad + \mathbb{1}_{[t<\tau]}\mathbb{1}_{[t\geq\tau_\alpha]}\frac{\sqrt{t-\bar{g}_t}}{\sqrt{T-\bar{g}_t}}\mathbb{E}[\mathbb{1}_{[\tau>T]}|\mathcal{G}_t]$$

$$(4.5) \qquad + \mathbb{1}_{[t<\tau_\alpha]}\frac{\alpha}{\sqrt{2}}\mathbb{E}\Big[\frac{1}{\sqrt{\tau-\bar{g}_{\tau_\alpha}}}\mathbb{1}_{[\tau\leq T]}\Big|\mathcal{G}_t\Big]$$

$$(4.6) \qquad + \mathbb{1}_{[t<\tau_\alpha]}\frac{\alpha}{\sqrt{2}}\mathbb{E}\Big[\frac{1}{\sqrt{T-\bar{g}_{\tau_\alpha}}}\mathbb{1}_{[\tau_\alpha\leq T<\tau]}\Big|\mathcal{G}_t\Big]$$

$$(4.7) \qquad + \mathbb{1}_{[t<\tau_\alpha]}\mathbb{E}[\mathbb{1}_{[\tau_\alpha>T]}|\mathcal{G}_t].$$

We next evaluate expressions (4.3)–(4.7). The distribution of the length of an excursion conditional on the age of the excursion is given in Chung (1976).

Conditional expectation in (4.3) on the event $[\tau_\alpha \leq t < \tau]$,

$$\mathbb{E}\left[\frac{1}{\sqrt{\tau - \bar{g}_t}} \mathbb{1}_{[\tau \leq T]} \Big| \mathcal{G}_t\right]$$

$$= \mathbb{E}\left[\frac{1}{\sqrt{L^\alpha}} \mathbb{1}_{[L^\alpha \leq T - \bar{g}_t]} \Big| \mathcal{G}_t\right]$$

$$= \int_{t - \bar{g}_t}^{T - \bar{g}_t} \frac{1}{2\sqrt{l}} \sqrt{\frac{t - \bar{g}_t}{l^3}} \, dl$$

$$= \frac{\sqrt{t - \bar{g}_t}}{2} \left(\frac{1}{t - \bar{g}_t} - \frac{1}{T - \bar{g}_t}\right).$$

Conditional expectation in (4.4) on the event $[\tau_\alpha \leq t < \tau]$,

$$\mathbb{E}[\mathbb{1}_{[\tau > T]} | \mathcal{G}_t] = \mathbb{Q}[T - \bar{g}_t < L^\alpha | \mathcal{G}_t]$$

$$= \int_{T - \bar{g}_t}^{\infty} \frac{1}{2} \sqrt{\frac{t - \bar{g}_t}{l^3}} \, dl$$

$$= \sqrt{\frac{t - \bar{g}_t}{T - \bar{g}_t}}.$$

Conditional expectation in (4.5) on the event $[\tau_\alpha > t]$,

$$\mathbb{E}\left[\frac{1}{\sqrt{\tau - \bar{g}_{\tau_\alpha}}} \mathbb{1}_{[\tau \leq T]} \Big| \mathcal{G}_t\right] = \mathbb{E}\left[\mathbb{E}\left[\frac{1}{\sqrt{\tau - \bar{g}_{\tau_\alpha}}} \mathbb{1}_{[\tau \leq T]} \Big| \mathcal{G}_{\tau_\alpha}\right] \Big| \mathcal{G}_t\right]$$

$$= \mathbb{E}\left[\mathbb{1}_{[T \geq \tau_\alpha]} \left(\frac{1}{\alpha\sqrt{2}} - \frac{\alpha}{2\sqrt{2}} \frac{1}{T - \bar{g}_{\tau_\alpha}}\right) \Big| \mathcal{G}_t\right].$$

Conditional expectation in (4.6) on the event $[\tau_\alpha > t]$,

$$\mathbb{E}\left[\frac{1}{\sqrt{T - \bar{g}_{\tau_\alpha}}} \mathbb{1}_{[\tau_\alpha \leq T < \tau]} \Big| \mathcal{G}_t\right] = \mathbb{E}\left[\mathbb{E}\left[\frac{1}{\sqrt{T - \bar{g}_{\tau_\alpha}}} \mathbb{1}_{[\tau_\alpha \leq T < \tau]} \Big| \mathcal{G}_{\tau_\alpha}\right] \Big| \mathcal{G}_t\right]$$

$$= \frac{\alpha}{\sqrt{2}} \mathbb{E}\left[\frac{1}{T - \bar{g}_{\tau_\alpha}} \mathbb{1}_{[\tau_\alpha \leq T]} \Big| \mathcal{G}_t\right].$$

Now, it remains to calculate $\mathbb{E}[\Delta V_\tau | \mathcal{G}_t]$ on $[t < \tau]$. Observe that $V_\tau = \mathbb{1}_{[\tau < T]} = \mathbb{1}_{[\tau \leq T]}$ since $\mathbb{Q}[\tau = T] = 0$. Thus, $\Delta V_\tau = \mathbb{1}_{[\tau \leq T]} - V_{\tau-}$. Since $\tau_\alpha < \tau$, a.s.,

$$V_{\tau-} = \left(\frac{\tau - \bar{g}_{\tau_\alpha}}{2} \left(\frac{1}{\tau - \bar{g}_{\tau_\alpha}} - \frac{1}{T - \bar{g}_{\tau_\alpha}}\right) + \frac{\tau - \bar{g}_{\tau_\alpha}}{T - \bar{g}_{\tau_\alpha}}\right) \mathbb{1}_{[\tau \leq T]}$$

$$= \frac{1}{2} \left(1 + \frac{\tau - \bar{g}_{\tau_\alpha}}{T - \bar{g}_{\tau_\alpha}}\right) \mathbb{1}_{[\tau \leq T]}.$$

Then,

$$\mathbb{E}[\Delta V_\tau | \mathcal{G}_t] \mathbb{1}_{[t < \tau]} = \frac{1}{2} \mathbb{E}[\mathbb{1}_{[\tau \leq T]} | \mathcal{G}_t] \mathbb{1}_{[t < \tau]}$$

$$- \frac{1}{2}\mathbb{E}\left[\frac{L^{\alpha}}{T - \bar{g}_{\tau_{\alpha}}}\mathbb{1}_{[\tau \leq T]}\Big|\mathcal{G}_t\right]\mathbb{1}_{[t < \tau]}$$

$$= \frac{1}{2}\mathbb{E}[\mathbb{1}_{[\tau \leq T]}|\mathcal{G}_t]\mathbb{1}_{[t < \tau]}$$

$$(4.8) \qquad - \frac{1}{2}\mathbb{E}\left[\frac{L^{\alpha}}{T - \bar{g}_{\tau_{\alpha}}}\mathbb{1}_{[\tau \leq T]}\Big|\mathcal{G}_t\right]\mathbb{1}_{[\tau_{\alpha} \leq t < \tau]}$$

$$(4.9) \qquad - \frac{1}{2}\mathbb{E}\left[\frac{L^{\alpha}}{T - \bar{g}_{\tau_{\alpha}}}\mathbb{1}_{[\tau \leq T]}\Big|\mathcal{G}_t\right]\mathbb{1}_{[t < \tau_{\alpha}]}.$$

Conditional expectation in (4.9) on the event $[t < \tau_{\alpha}]$,

$$\mathbb{E}\left[\frac{L^{\alpha}}{T - \bar{g}_{\tau_{\alpha}}}\mathbb{1}_{[\tau \leq T]}\Big|\mathcal{G}_t\right] = \mathbb{E}\left[\mathbb{E}\left[\frac{L^{\alpha}}{T - \bar{g}_{\tau_{\alpha}}}\mathbb{1}_{[\tau \leq T]}\Big|\mathcal{G}_{\tau_{\alpha}}\right]\Big|\mathcal{G}_t\right]$$

$$(4.10) \qquad\qquad = \mathbb{E}\left[\left(\frac{\int_{\alpha^2/2}^{T - \bar{g}_{\tau_{\alpha}}} l/2\sqrt{(\alpha^2/2)l^{-3}}\,dl}{T - \bar{g}_{\tau_{\alpha}}}\right)\mathbb{1}_{[\tau_{\alpha} \leq T]}\Big|\mathcal{G}_t\right]$$

$$= \mathbb{E}\left[\frac{\alpha/\sqrt{2}(\sqrt{T - \bar{g}_{\tau_{\alpha}}} - \alpha/\sqrt{2})}{T - \bar{g}_{\tau_{\alpha}}}\mathbb{1}_{[\tau_{\alpha} \leq T]}\Big|\mathcal{G}_t\right].$$

Similarly, conditional expectation in (4.8) on the event $[\tau_{\alpha} \leq t < \tau]$,

$$\mathbb{E}\left[\frac{L^{\alpha}}{T - \bar{g}_{\tau_{\alpha}}}\mathbb{1}_{[\tau \leq T]}\Big|\mathcal{G}_t\right] = \frac{\alpha/\sqrt{2}(\sqrt{T - \bar{g}_t} - \alpha/\sqrt{2})}{T - \bar{g}_t}.$$

Therefore,

$$\mathbb{E}[\mathbb{1}_{[\tau > T]}|\mathcal{G}_t]\mathbb{1}_{[t < \tau]}$$

$$= \mathbb{1}_{[t < \tau]}\mathbb{1}_{[t \geq \tau_{\alpha}]}\frac{1}{2}\left(1 - \frac{t - \bar{g}_t}{T - \bar{g}_t}\right)$$

$$+ \mathbb{1}_{[t < \tau]}\mathbb{1}_{[t \geq \tau_{\alpha}]}\frac{t - \bar{g}_t}{T - \bar{g}_t}$$

$$+ \mathbb{1}_{[t < \tau_{\alpha}]}\frac{\alpha}{\sqrt{2}}\mathbb{E}\left[\mathbb{1}_{[T \geq \tau_{\alpha}]}\left(\frac{1}{\alpha\sqrt{2}} - \frac{\alpha}{2\sqrt{2}}\frac{1}{T - \bar{g}_{\tau_{\alpha}}}\right)\Big|\mathcal{G}_t\right]$$

$$(4.11) \qquad + \mathbb{1}_{[t < \tau_{\alpha}]}\frac{\alpha^2}{2}\mathbb{E}\left[\frac{1}{T - \bar{g}_{\tau_{\alpha}}}\mathbb{1}_{[\tau_{\alpha} \leq T]}\Big|\mathcal{G}_t\right]$$

$$+ \mathbb{1}_{[t < \tau_{\alpha}]}\mathbb{E}[\mathbb{1}_{[\tau_{\alpha} > T]}|\mathcal{G}_t]$$

$$- \frac{1}{2}\mathbb{1}_{[t < \tau]}\mathbb{E}[\mathbb{1}_{[\tau \leq T]}|\mathcal{G}_t]$$

$$+ \mathbb{1}_{[\tau_{\alpha} \leq t < \tau]}\frac{\alpha/(2\sqrt{2})(\sqrt{T - \bar{g}_t} - \alpha/\sqrt{2})}{T - \bar{g}_t}$$

$$+ \mathbb{1}_{[t < \tau_\alpha]} \mathbb{E}\left[ \frac{\alpha/(2\sqrt{2})(\sqrt{T - \bar{g}_{\tau_\alpha}} - \alpha/\sqrt{2})}{T - \bar{g}_{\tau_\alpha}} \mathbb{1}_{[\tau_\alpha \le T]} \Big| \mathcal{G}_t \right],$$

which yields

$$
\begin{aligned}
\mathbb{E}[\mathbb{1}_{[\tau > T]} | \mathcal{G}_t] \mathbb{1}_{[t < \tau]} = & -\mathbb{1}_{[t < \tau]} + \mathbb{1}_{[t < \tau]} \mathbb{1}_{[t \ge \tau_\alpha]} \left( 1 + \frac{t - \bar{g}_t}{T - \bar{g}_t} \right) \\
& + \mathbb{1}_{[t < \tau_\alpha]} \mathbb{E}\left[ \mathbb{1}_{[T \ge \tau_\alpha]} \left( 1 + \frac{\alpha}{\sqrt{2}} \frac{1}{\sqrt{T - \bar{g}_{\tau_\alpha}}} \right) \Big| \mathcal{G}_t \right] \\
& + \mathbb{1}_{[t < \tau_\alpha]} 2 \mathbb{E}[\mathbb{1}_{[\tau_\alpha > T]} | \mathcal{G}_t] \\
& + \mathbb{1}_{[\tau_\alpha \le t < \tau]} \frac{\alpha/\sqrt{2}(\sqrt{T - \bar{g}_t} - \alpha/\sqrt{2})}{T - \bar{g}_t}.
\end{aligned}
$$

$(4.12)$

In order to get the price, $S_t$, we need to obtain the law of $\tau_\alpha$ on the event $[t < \tau_\alpha]$ conditional on $\mathcal{G}_t$. To find the Laplace transform of this density we introduce the following martingale as in Chesney, Jeanblanc-Picqué and Yor ([1997](#)):

$$N_t := \Psi(-\lambda \mu_{t \wedge \tau_\alpha}) \exp\left( -\frac{\lambda^2}{2}(t \wedge \tau_\alpha) \right),$$

where $\mu_t = M_t/\sqrt{2}$ and $\Psi(z) = \int_0^\infty x \exp(zx - x^2/2)\, dx$. Using the optional stopping theorem, we obtain

$$
\begin{aligned}
\mathbb{E}\left[ \Psi(-\lambda \mu_{\tau_\alpha}) \exp\left( -\frac{\lambda^2}{2}(\tau_\alpha) \right) \Big| \mathcal{G}_t \right] \\
= \Psi(-\lambda \mu_{t \wedge \tau_\alpha}) \exp\left( -\frac{\lambda^2}{2}(t \wedge \tau_\alpha) \right),
\end{aligned}
$$

which in turn implies

$$\mathbb{1}_{[t > \tau_\alpha]} \mathbb{E}\left[ \exp\left( -\frac{\lambda^2}{2} \tau_\alpha \right) \Big| \mathcal{G}_t \right] = \mathbb{1}_{[t > \tau_\alpha]} \frac{\Psi(-\lambda \mu_t) \exp((-\lambda^2/2)t)}{\Psi(\lambda \alpha/\sqrt{2})}.$$

Using relatively standard software, one can invert this Laplace transform and compute the expectations given in expressions [(4.3)](#)–[(4.12)](#).

For time 0, using expression [(4.12)](#) and rearranging the terms give

$$(4.13) \quad S_0 = \exp\left( -\int_0^T r_u\, du \right) \left( 1 - \left( \mathbb{Q}[\tau_\alpha \le T] - \mathbb{E}\left[ \frac{\alpha/\sqrt{2}}{\sqrt{T - \bar{g}_{\tau_\alpha}}} \mathbb{1}_{[\tau_\alpha \le T]} \right] \right) \right).$$

This is the price of the risky zero coupon bond at time 0. The interpretation of the last term in this expression is important. Default occurs not at time $\tau_\alpha$, but at time $\tau$. The default time $\tau$ is, therefore, less likely than the hitting time $\tau_\alpha$. The probability $\mathbb{Q}[\tau_\alpha \le T]$ is reduced to account for this difference.

Unfortunately, the law of $\tau_\alpha$ is only known through its Laplace transform, which is very difficult to invert analytically. See Chesney, Jeanblanc-Picqué and Yor (1997) in this respect, which gives the following formula:

$$\mathbb{E}\left[\exp\left(-\frac{\lambda^2}{2}\tau_\alpha\right)\right] = \frac{1}{\Psi(\lambda\alpha/\sqrt{2})}.$$

Inverting this Laplace transform yields the law for $\tau_\alpha$, and given the law for $\tau_\alpha$, expression (4.13) is easily computed.

**5. Conclusion.** This paper provides an alternative method for generating reduced form credit risk models from structural models. The difference from Duffie and Lando (2001) is that instead of using filtering theory to go from the manager's information to the market's as in Duffie and Lando, we use a reduction of the manager's information set. This modification is both conceptually and mathematically a different approach to the topic. Indeed, the perspective from filtering theory is that the market's information set is the same as the manager's, but with additional noise included. The perspective from reducing the manager's information set is that the market's information set is the same as the manager's, but the market just knows less of it. It would be interesting to investigate more complex structural models than those used herein and more complex information reductions.

**Acknowledgment.** The authors thank Professor Monique Jeanblanc-Picqué for pointing out errors in earlier versions of this paper.

U. Cetin
Center for Applied Mathematics
Cornell University
Ithaca, New York 14853-3801
USA
e-mail: umut@cam.cornell.edu

P. Protter
Operations Research
    and Industrial Engineering Department
Cornell University
Ithaca, New York 14853-3801
USA
e-mail: protter@orie.cornell.edu

R. Jarrow
Johnson Graduate School
    of Management
Cornell University
Ithaca, New York 14853-3801
USA
e-mail: raj15@cornell.edu

Y. Yildirim
School of Management
Syracuse University
Syracuse, New York 13244
USA
e-mail: yildiray@syr.edu